\documentclass[12pt]{article}
\usepackage{titlesec}
\titleformat*{\section}{\Large\bfseries\fontfamily{qhv}\selectfont}
\titleformat*{\subsection}{\large\bfseries\fontfamily{qhv}\selectfont}
\usepackage{graphicx}
\usepackage{hyperref}
\usepackage{tikz}
\usetikzlibrary{calc}
\usetikzlibrary{positioning}
\usepackage[margin=0.75in]{geometry}
\usepackage{mathrsfs}
\usepackage{amsfonts}
\usepackage{amsmath}
\usepackage{amssymb}
\usepackage{amsthm}
\usepackage{mdframed}
\usepackage{enumitem}
\usepackage{wallpaper}
\usepackage{hyperref}
\usepackage[misc]{ifsym}
\usepackage{float}
\usepackage[hang,flushmargin]{footmisc}

\newcommand{\T}{\mathcal{T}}
\newcommand{\V}{\mathcal{V}}
\newcommand{\E}{\mathcal{E}}
\newcommand{\ds}{\displaystyle}

\usepackage{subcaption}
\usepackage{graphicx}

\DeclareCaptionFormat{custom}
{%
    \textbf{\fontfamily{qhv}\selectfont{#1 }}{ #3}
}
\newcommand\blfootnote[1]{%
  \begingroup
  \renewcommand\thefootnote{}\footnote{#1}%
  \addtocounter{footnote}{-1}%
  \endgroup
}

\makeatletter
\renewenvironment{proof}[1][{\fontfamily{qhv}\selectfont{\proofname}}] {\begin{mdframed}\par\pushQED{\qed}\normalfont\topsep6\p@\@plus6\p@\relax\trivlist\item[\hskip\labelsep\bfseries#1\@addpunct{}]$\newline$\ignorespaces}{\popQED\endtrivlist\@endpefalse\end{mdframed}}
\makeatother

\newtheoremstyle{colon}%
{}
{}
{}
{}
{\bfseries\fontfamily{qhv}\selectfont}
{\phantom{:}}
{ }
{}

\theoremstyle{colon}
\newtheorem{thm}{Theorem}

\newtheorem{defn}{Definition}

\allowdisplaybreaks

\begin{document}
\title{\textbf{{\fontfamily{qhv}\selectfont{Packing and Covering Triangles in Bilaterally-
Complete Tripartite Graphs}}}}
	\author{Naivedya Amarnani $\pmb{\cdot}$ Amaury De
Burgos $\pmb{\cdot}$ Wayne Broughton}
  \maketitle
  
 \begin{abstract}
 	We use Menger's Theorem and K\"onig's Line Colouring Theorem to show that in any tripartite graph with two complete (bipartite) sides the maximum number of pairwise edge-disjoint triangles equals the minimum number of edges that meet all triangles. This generalizes the corresponding result for complete tripartite graphs given by Lakshmanan, et al.
    
  \vspace{6pt}
  \noindent\textbf{\fontfamily{qhv}\selectfont{Keywords:}} Maximum Packing, Tripartite Graphs, Minimum T-Transversal, Triangles, Complete Bipartite Graph
 \end{abstract}
 
 \section{Introduction}
\label{intro} \blfootnote{The first author was supported by Mitacs through a Mitacs Globalink Research Internship while visiting UBC Okanagan. The second author was supported by an Undergraduate Research Award from the I.K. Barber Faculty of Science at UBC Okanagan.\vspace{-2mm}\\\rule{7.08cm}{0.4pt}\\N. Amarnani\\
Mathematics and Statistics, Indian Institute of Technology Kanpur\\
Kanpur, Uttar Pradesh, India\\\\A. De Burgos, W. Broughton\\
Mathematics, University of British Columbia - Okanagan\\
Kelowna, British Columbia, Canada\\
E-mail: wayne.broughton@ubc.ca
}
For standard graph-theoretic terms or notation, refer to \cite{DW}. All graphs will be finite, simple and undirected unless indicated otherwise. An edge with endpoints $u$ and $v$ will generally be denoted $uv$. In graph $G=(V,E)$, if $X\subseteq V$ or $Y\subseteq E$ then $G[X]$ denotes the subgraph of $G$ induced by $X$ and $G[Y]$ denotes the subgraph edge-induced by $Y$.

A \textit{triangle} in a graph $G$ is a subgraph isomorphic to $K_3$, consisting of 3 vertices that are all joined by edges. Most of the following definitions and notation are from \cite{YR} and \cite{LBT}.

Let $\T$ denote the set of all triangles in $G=(V,E)$. A subset $E'\subseteq E$ is said to be a \textit{$\T$-transversal} if every triangle in $\T$ has an edge in $E'$. 	

A subset $\T'\subseteq \T$ is said to be a \textit{packing} in $G$ if all the triangles in $\T'$ are pairwise edge-disjoint. 

We denote the cardinality of the smallest $\T$-transversal in $G$ by $\tau_{\triangle}(G)$ and the cardinality of the largest packing by $\nu_{\triangle}(G)$. It is not difficult to see that $\nu_{\triangle}(G)\leq\tau_{\triangle}(G)\leq3\nu_{\triangle}(G)$.

In 1981 Zsolt Tuza conjectured (\cite{T}) that $\tau_{\triangle}(G)\leq2\nu_{\triangle}(G)$ for every graph $G$. The general case is still open, but when $G$ is a tripartite graph an even stronger bound was proved by Haxell and Kohayakawa (\cite{HK}). In 2012, Lakshmanan et al. (\cite{LBT}) proved that in fact the equality $\tau_{\triangle}(G)=\nu_{\triangle}(G)$ holds for certain classes of tripartite graphs. Their results can be viewed as extensions to these tripartite graphs of the K\"onig-Egerv\'ary Theorem, equating the minimum size of a vertex cover and the maximum matching size in bipartite graphs. In particular, it is observed in \cite{LBT} that $\tau_{\triangle}(G)=\nu_{\triangle}(G)$ always holds for complete tripartite graphs. In this paper, we show that this equality also holds for tripartite graphs with only two complete bipartite sides.

If $G=(V,E)$ is a tripartite graph whose vertex set is partitioned into  $V=A\cup B\cup C$, then we write $G=(A,B,C;E)$ and we call the induced bipartite subgraphs $G[A\cup B]$, $G[A\cup C]$, and $G[B\cup C]$ the \textit{sides} of $G$. The sets of edges on the sides are written $E_{AB}$, $E_{AC}$, and $E_{BC}$, respectively. More generally, if $M,W\subseteq V$ then $E_{MW}\subseteq E$ denotes the set of edges with one endpoint in each of $M$ and $W$ (abusing notation slightly if $M$ or $W$ is just a single vertex).

\begin{defn}[Bilaterally-Complete Tripartite Graph]
Let $G=(A,B,C;E)$ be a tripartite graph. If (at least) two of the sides of $G$ are complete bipartite graphs, then we call $G$ a \textit{bilaterally-complete tripartite graph}.	
\end{defn}

Our main result is the following:
\begin{thm}\label{thm:1}
\textit{If $G$ is a bilaterally-complete tripartite graph, then $\tau_{\triangle}(G)=\nu_{\triangle}(G)$.}
\end{thm}

\textbf{\fontfamily{qhv}\selectfont{Note:}} It is possible to derive a proof of this fact from a complicated chain of results in some old papers (particularly \cite{CM}, \cite{DS} and \cite{HP}), but we give a straightforward proof relying only on classical theorems. The alternative approach will be discussed in the last section of this paper.

\section{Characterization of Minimum $\pmb{\T}$-transversals}
\label{SMTT}
\hypertarget{SMTT}{}
Let $G=(A,B,C;E)$ be a bilaterally-complete tripartite graph with complete sides $G[A\cup B]\cong K_{p,q}$ and $G[A\cup C]\cong K_{p,r}$, where $|A|=p$, $|B|=q$, $|C|=r$. We will describe the structure of a minimum $\T$-transversal of $G$ in terms of vertices and edges of $G[B\cup C]$.
	
	For each vertex $a\in A$, the triangles in $G$ containing $a$ are in one-to-one correspondence with the edges in $E_{BC}$, since $G$ is bilaterally-complete. Likewise, edges in $E_{aB}$ and $E_{aC}$ can be identified with vertices in $B$ and $C$, respectively. 
	
	Now let $E'$ be a $\T$-transversal in $G$, and fix $a\in A$. If we let $E'_{BC}=E'\cap E_{BC}$, then the triangles containing $a$ and corresponding to edges in $E_{BC}\setminus E'_{BC}$ must be covered by $E'$ using edges in $E_{aB}\cup E_{aC}$. This means the corresponding vertices in $B\cup C$ must cover all the edges in $E_{BC}\setminus E'_{BC}$; that is, they form a vertex cover $W_a$ of $G[E_{BC}\setminus E'_{BC}]$ for each $a\in A$. 	
 
Conversely, if we choose a set of edges $E'_{BC}\subseteq E_{BC}$ and a vertex cover $W_a$ of $G[E_{BC}\setminus E'_{BC}]$ for each $a\in A$, then 
	\[
	E'=E'_{BC}\cup  \bigcup_{a\in A} E_{aW_a}
	\] 
will be a $\T$-transversal of $G$ with size $|E'|=|E'_{BC}|+\sum_{a\in A}|W_a|$.

	If $W_a=W$ is the same vertex cover for all $a\in A$ we will call this a \emph{uniform} $\T$-transversal.
	Without loss of generality, a minimum $\T$-transversal can be assumed to be uniform by choosing $E'_{BC}\subseteq E_{BC}$ and a minimum vertex cover $W$ of $G[E_{BC}\setminus E'_{BC}]$, such that the value of $|E'_{BC}|+p|W|$ is minimized.
	Hence,
	\begin{align*}
		\tau_{\triangle}(G)=\min_{E'_{BC}\subseteq E_{BC}}\{|E'_{BC}|+p|W|\},
	\end{align*}
	where $W$ is a minimum vertex cover of $G[E_{BC}\setminus E'_{BC}]$.
	
	It is worth pointing out that in bilaterally-complete tripartite graphs, the edges of the smallest side do not always form a minimum $\mathcal{T}$-transversal. In fact, there exist bilaterally-complete tripartite graphs whose minimum $\mathcal{T}$-transversals must have an edge in every side of the graph. Such a graph is pictured in \hypertarget{G}{Figure 1}. This makes the bilaterally-complete case more complicated than the complete tripartite case, where the smallest side is always a minimum $\mathcal{T}$-transversal and the maximum packing can be described explicitly as a selection of one triangle on each edge in the smallest side of the graph.	

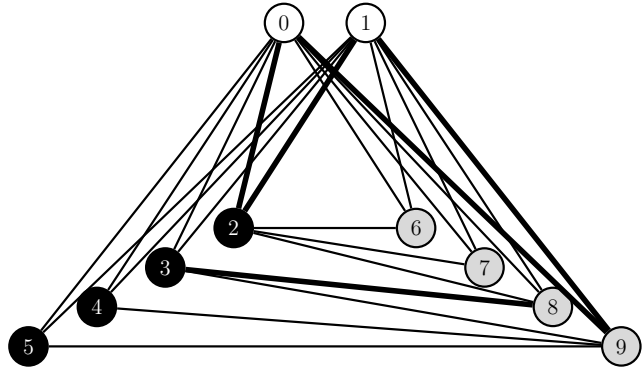
\begin{figure}[H]
\hypertarget{G}{}
{\centering\scalebox{0.7}{%
			\begin{tikzpicture}
			[
				W/.style={circle, draw=black, fill=none, very thick, minimum size=7mm},
				B/.style={circle, draw=black, fill=black, very thick, minimum size=7mm, text=white},
				G/.style={circle, draw=black, fill=black!15, very thick, minimum size=7mm},
				S/.style={circle, draw=none, fill=none, very thick, minimum size=7mm}
			]
			\node(o){};
			\node[W](0) at ([shift=({105:3cm})]o) {$0$};
			\node[W](1) at ([shift=({75:3cm})]o) {$1$};
			\node[B](2) at ([shift=({210:2cm})]o) {$2$};
			\node[B](3) at ([shift=({210:3.5cm})]o) {$3$};
			\node[B](4) at ([shift=({210:5cm})]o) {$4$};
			\node[B](5) at ([shift=({210:6.5cm})]o) {$5$};
			\node[G](6) at ([shift=({330:2cm})]o) {$6$};
			\node[G](7) at ([shift=({330:3.5cm})]o) {$7$};
			\node[G](8) at ([shift=({330:5cm})]o) {$8$};
			\node[G](9) at ([shift=({330:6.5cm})]o) {$9$};
			\draw[line width = 3pt](0)--(2);
			\draw[very thick](0)--(3);
			\draw[very thick](0)--(4);
			\draw[very thick](0)--(5);
			
			\draw[line width = 3pt](1)--(2);
			\draw[very thick](1)--(3);
			\draw[very thick](1)--(4);
			\draw[very thick](1)--(5);
			
			\draw[line width = 3pt](0)--(9);
			\draw[very thick](0)--(7);
			\draw[very thick](0)--(8);
			\draw[very thick](0)--(6);
			
			\draw[line width = 3pt](1)--(9);
			\draw[very thick](1)--(7);
			\draw[very thick](1)--(8);
			\draw[very thick](1)--(6);
			
			\draw[very thick](2)--(6);
			\draw[very thick](2)--(7);
			\draw[very thick](2)--(8);
			\draw[line width = 3pt](3)--(8);
			\draw[very thick](3)--(9);
			\draw[very thick](4)--(9);
			\draw[very thick](5)--(9);
	\end{tikzpicture}}%
	
}
\caption{A bilaterally-complete tripartite graph $G=(A,B,C;E)$ with $G[A\cup B]\cong G[A\cup C]\cong K_{2,4}$. The smallest side has seven edges, but the only minimum $\mathcal{T}$-transversal is the one of size five shown in boldface, with an edge in every side of the graph. In this example $E'_{BC}$ is just the singleton $\{38\}$ and $W=\{2,9\}$.}
\label{fig:1}
\end{figure}

\section{Network Graph \textit{H(G)}}
\label{NWG}
As before, let $G=(A,B,C;E)$ be a bilaterally-complete tripartite graph with complete sides $G[A\cup B]\cong K_{p,q}$ and $G[A\cup C]\cong K_{p,r}$, where $|A|=p$, $|B|=q$, $|C|=r$. The \textit{network graph} $H(G)=(\V,\E)$ is a directed graph constructed from $G[B\cup C]$ as follows:
\begin{enumerate}
			\item Make $p$ copies of each of the vertex sets $B$ and $C$, labelled $B_1,\ldots,B_p$ and $C_1,\ldots,C_p$.
			\item Create a vertex for each edge in $E_{BC}$, and denote this set of vertices by $V_{BC}$. 
			\item Create a source vertex $s$ and a sink vertex $t$. 
			\item For each edge $uv$ in $E_{BC}$ ($u\in B$), make a directed arc from every copy of $u$ in  $B_1,\ldots,B_p$ to the vertex in $V_{BC}$ corresponding to the edge $uv$, and an arc from that same vertex in $V_{BC}$ to every copy of $v$ in $C_1,\ldots,C_p$.
			\item Make a directed arc from $s$ to each vertex in $B_1\cup\cdots\cup B_p$, and from every vertex in $C_1\cup\cdots\cup C_p$ to $t$.
\end{enumerate}

As an example, let $G$ be the bilaterally-complete tripartite graph shown in \hyperlink{G}{Figure 1}. Its induced bipartite graph $G[B\cup C]$ is shown in \hyperlink{BG}{Figure 2} and the resulting network graph $H(G)$ is shown in \hyperlink{HG}{Figure 3}. The general representation of $H(G)$ is shown in \hyperlink{GHG}{Figure 4}.

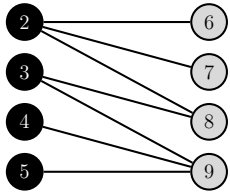
\begin{figure}[H]
\hypertarget{BG}{}
{\centering\scalebox{0.65}{%
			\begin{tikzpicture}
			[
				W/.style={circle, draw=black, fill=none, very thick, minimum size=7mm},
				B/.style={circle, draw=black, fill=black, very thick, minimum size=7mm, text=white},
				G/.style={circle, draw=black, fill=black!15, very thick, minimum size=7mm},
				S/.style={circle, draw=none, fill=none, very thick, minimum size=7mm}
			]
			\node[B](2){$2$};
			\node[B](3)[below =0.25cm of 2]{$3$};
			\node[B](4)[below =0.25cm of 3]{$4$};
			\node[B](5)[below =0.25cm of 4]{$5$};
			
			\node[G](6)[right =3cm of 2]{$6$};
			\node[G](7)[right =3cm of 3]{$7$};
			\node[G](8)[right =3cm of 4]{$8$};
			\node[G](9)[right =3cm of 5]{$9$};
			\draw[very thick](2)--(6);
			\draw[very thick](2)--(7);
			\draw[very thick](2)--(8);
			\draw[very thick](3)--(8);
			\draw[very thick](3)--(9);
			\draw[very thick](4)--(9);
			\draw[very thick](5)--(9);
	\end{tikzpicture}}%

}
\caption{The induced bipartite graph $G[B\cup C]$ of the graph in Figure~\ref{fig:1}.}
\label{fig:2}
\end{figure}

\begin{figure*}[ht]
\hypertarget{HG}{}
\begin{center}\scalebox{0.65}{%
			\begin{tikzpicture}%
			[
				roundnode/.style={circle, draw=black, fill=none, very thick, minimum size=7mm, inner sep = 0cm},
				W/.style={circle, draw=black, fill=none, very thick, minimum size=7mm},
				B/.style={circle, draw=black, fill=black, very thick, minimum size=7mm, text=white},
				G/.style={circle, draw=black, fill=black!15, very thick, minimum size=7mm},
				S/.style={circle, draw=none, fill=none, very thick, minimum size=7mm}
			]
			\node[B](2){$2$};
			\node[B](3)[below =0.25cm of 2]{$3$};
			\node[B](4)[below =0.25cm of 3]{$4$};
			\node[B](5)[below =0.25cm of 4]{$5$};
			\node[S](s)[below =0.5 of 5]{\phantom{s}};
			\node[B](2')[below =0.5 cm of s]{$2'$};
			\node[B](3')[below =0.25cm of 2']{$3'$};
			\node[B](4')[below =0.25cm of 3']{$4'$};
			\node[B](5')[below =0.25cm of 4']{$5'$};
			
			\node[G](6)[right =6cm of 2]{$6$};
			\node[G](7)[right =6cm of 3]{$7$};
			\node[G](8)[right =6cm of 4]{$8$};
			\node[G](9)[right =6cm of 5]{$9$};
			\node[S](t)[right =6cm of s]{\phantom{t}};
			\node[G](6')[right =6cm of 2']{$6'$};
			\node[G](7')[right =6cm of 3']{$7'$};
			\node[G](8')[right =6cm of 4']{$8'$};
			\node[G](9')[right =6cm of 5']{$9'$};
			
			\node[roundnode](38) at ($(s)!0.5!(t)$) {$38$};
			\node[roundnode](28) [above =0.25 cm of 38]{$28$};
			\node[roundnode](27) [above =0.25 cm of 28]{$27$};
			\node[roundnode](26) [above =0.25 cm of 27]{$26$};
			\node[roundnode](39) [below =0.25 cm of 38]{$39$};
			\node[roundnode](49) [below =0.25 cm of 39]{$49$};
			\node[roundnode](59) [below =0.25 cm of 49]{$59$};
			
			\node[W](S)[left =3cm of s]{\normalsize $s$};
			\node[W](T)[right =3cm of t]{\normalsize $t$};
			\draw[very thick,-latex](2)--(26);
			\draw[very thick,-latex](26)--(6);
			\draw[very thick,-latex](2')--(26);
			\draw[very thick,-latex](26)--(6');
			\draw[very thick,-latex](2)--(27);
			\draw[very thick,-latex](27)--(7);
			\draw[very thick,-latex](2')--(27);
			\draw[very thick,-latex](27)--(7');
			\draw[very thick,-latex](2)--(28);
			\draw[very thick,-latex](28)--(8);
			\draw[very thick,-latex](2')--(28);
			\draw[very thick,-latex](28)--(8');
			\draw[very thick,-latex](3)--(38);
			\draw[very thick,-latex](38)--(8);
			\draw[very thick,-latex](3')--(38);
			\draw[very thick,-latex](38)--(8');
			\draw[very thick,-latex](3)--(39);
			\draw[very thick,-latex](39)--(9);
			
			\draw[very thick,-latex](3')--(39);
			\draw[very thick,-latex](4)--(49);
			\draw[very thick,-latex](4')--(49);
			\draw[very thick,-latex](5)--(59);
			\draw[very thick,-latex](5')--(59);

			\draw[very thick,-latex](39)--(9');
			\draw[very thick,-latex](49)--(9);
			\draw[very thick,-latex](49)--(9');
			\draw[very thick,-latex](59)--(9);
			\draw[very thick,-latex](59)--(9');

			\draw[very thick,-latex](S)--(2);
			\draw[very thick,-latex](S)--(3);
			\draw[very thick,-latex](S)--(4);
			\draw[very thick,-latex](S)--(5);
			\draw[very thick,-latex](S)--(2');
			\draw[very thick,-latex](S)--(3');
			\draw[very thick,-latex](S)--(4');
			\draw[very thick,-latex](S)--(5');
			
			\draw[very thick,-latex](6)--(T);
			\draw[very thick,-latex](7)--(T);
			\draw[very thick,-latex](8)--(T);
			\draw[very thick,-latex](9)--(T);
			\draw[very thick,-latex](6')--(T);
			\draw[very thick,-latex](7')--(T);
			\draw[very thick,-latex](8')--(T);
			\draw[very thick,-latex](9')--(T);
	\end{tikzpicture}}
\end{center}
\caption{The network graph $H(G)$ constructed from $G[B\cup C]$ in Figure~\ref{fig:2}.}
\label{fig:3}
\end{figure*}
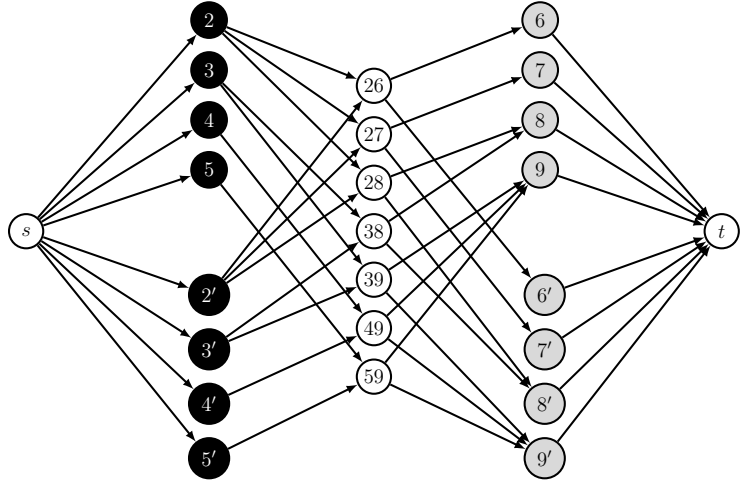

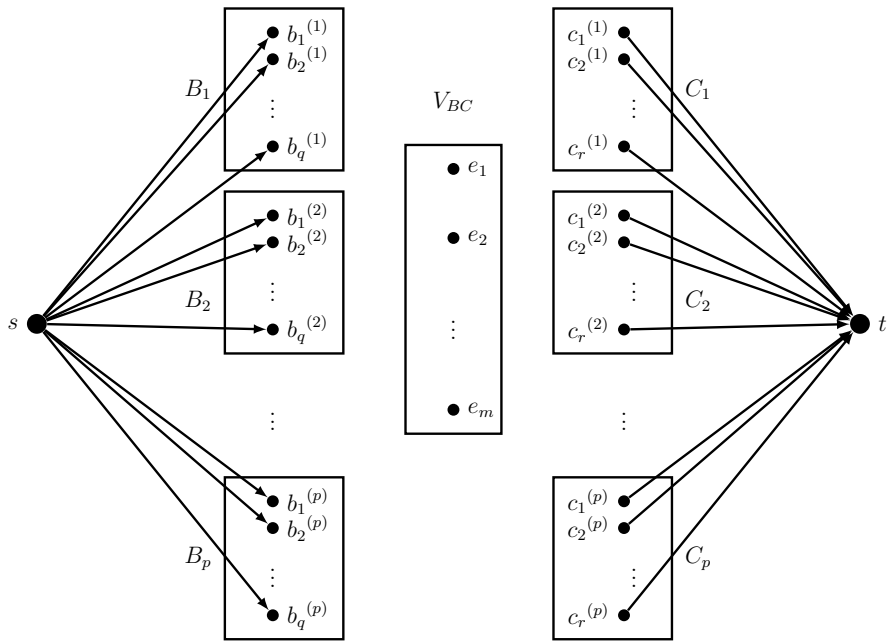
\begin{figure*}[h!]
\hypertarget{GHG}{}
\begin{center}\scalebox{0.75}{%
			\begin{tikzpicture}
			\node(b1)[circle,fill,inner sep=0pt,minimum size=6pt, label = right:{${b_1}^{(1)}$}]{};
			\node(b2)[circle,fill,inner sep=0pt,minimum size=6pt, label = right:{${b_2}^{(1)}$}][below =0.25cm of b1]{};
			\node(bdots)[below =0.25cm of b2]{$\vdots$};
			\node(bq)[circle,fill,inner sep=0pt,minimum size=6pt, label = right:{${b_q}^{(1)}$}][below =0.25cm of bdots]{};
			\node(b1')[circle,fill,inner sep=0pt,minimum size=6pt, label = right:{${b_1}^{(2)}$}][below = of bq]{};
			\node(b2')[circle,fill,inner sep=0pt,minimum size=6pt, label = right:{${b_2}^{(2)}$}][below =0.25cm of b1']{};
			\node(bdots')[below =0.25cm of b2']{$\vdots$};
			\node(bq')[circle,fill,inner sep=0pt,minimum size=6pt, label = right:{${b_q}^{(2)}$}][below =0.25cm of bdots']{};
			\node(BDOTS)[below = of bq']{$\vdots$};
			\node(b1'')[circle,fill,inner sep=0pt,minimum size=6pt, label = right:{${b_1}^{(p)}$}][below = of BDOTS]{};
			\node(b2'')[circle,fill,inner sep=0pt,minimum size=6pt, label = right:{${b_2}^{(p)}$}][below =0.25cm of b1'']{};
			\node(bdots'')[below =0.25cm of b2'']{$\vdots$};
			\node(bq'')[circle,fill,inner sep=0pt,minimum size=6pt, label = right:{${b_q}^{(p)}$}][below =0.25cm of bdots'']{};
			\draw[very thick, -]($(b1)+(-0.85,0.425)$) rectangle node[left = 1.125cm]{$B_1$} ($(bq)+(1.25,-0.425)$);
			\draw[very thick, -]($(b1')+(-0.85,0.425)$) rectangle node[left = 1.125cm,yshift=-0.5cm]{$B_2$}($(bq')+(1.25,-0.425)$);
			\draw[very thick, -]($(b1'')+(-0.85,0.425)$) rectangle node[left = 1.125cm]{$B_p$}($(bq'')+(1.25,-0.425)$);
			\node(s)[circle,fill,inner sep=0pt,minimum size=10pt, label = left:{\normalsize $s$}][left =4cm of $(b1)!0.5!(bq'')$]{};
			\draw[very thick,-latex](s)--(b1);
			\draw[very thick, -latex](s)--(b2);
			\draw[very thick, -latex](s)--(bq);
			\draw[very thick, -latex](s)--(b1');
			\draw[very thick, -latex](s)--(b2');
			\draw[very thick, -latex](s)--(bq');
			\draw[very thick, -latex](s)--(b1'');
			\draw[very thick, -latex](s)--(b2'');
			\draw[very thick, -latex](s)--(bq'');
			\node(c1)[circle,fill,inner sep=0pt,minimum size=6pt, label = left:{${c_1}^{(1)}$}][right =6cm of b1]{};
			\node(c2)[circle,fill,inner sep=0pt,minimum size=6pt, label = left:{${c_2}^{(1)}$}][right =6cm of b2]{};
			\node(cdots)[right =6cm of bdots]{$\vdots$};
			\node(cr)[circle,fill,inner sep=0pt,minimum size=6pt, label = left:{${c_r}^{(1)}$}][right =6cm of bq]{};
			\node(c1')[circle,fill,inner sep=0pt,minimum size=6pt, label = left:{${c_1}^{(2)}$}][right =6cm of b1']{};
			\node(c2')[circle,fill,inner sep=0pt,minimum size=6pt, label = left:{${c_2}^{(2)}$}][right =6cm of b2']{};
			\node(cdots')[right =6cm of bdots']{$\vdots$};
			\node(cr')[circle,fill,inner sep=0pt,minimum size=6pt, label = left:{${c_r}^{(2)}$}][right =6cm of bq']{};
			\node(CDOTS)[below = of cr']{$\vdots$};
			\node(c1'')[circle,fill,inner sep=0pt,minimum size=6pt, label = left:{${c_1}^{(p)}$}][right =6cm of b1'']{};
			\node(c2'')[circle,fill,inner sep=0pt,minimum size=6pt, label = left:{${c_2}^{(p)}$}][right =6cm of b2'']{};
			\node(cdots'')[right =6cm of bdots'']{$\vdots$};
			\node(cr'')[circle,fill,inner sep=0pt,minimum size=6pt, label = left:{${c_r}^{(p)}$}][right =6cm of bq'']{};
			\draw[very thick, -]($(c1)+(-1.25,0.425)$) rectangle node[right = 1.125cm]{$C_1$} ($(cr)+(0.85,-0.425)$);
			\draw[very thick, -]($(c1')+(-1.25,0.425)$) rectangle node[right = 1.125cm,yshift=-0.5cm]{$C_2$}($(cr')+(0.85,-0.425)$);
			\draw[very thick, -]($(c1'')+(-1.25,0.425)$) rectangle node[right = 1.125cm]{$C_p$}($(cr'')+(0.85,-0.425)$);
			\node(t)[circle,fill,inner sep=0pt,minimum size=10pt, label = right:{\normalsize $t$}][right =4cm of $(c1)!0.5!(cr'')$]{};
			\draw[very thick,latex-](t)--(c1);
			\draw[very thick, latex-](t)--(c2);
			\draw[very thick, latex-](t)--(cr);
			\draw[very thick, latex-](t)--(c1');
			\draw[very thick, latex-](t)--(c2');
			\draw[very thick, latex-](t)--(cr');
			\draw[very thick, latex-](t)--(c1'');
			\draw[very thick, latex-](t)--(c2'');
			\draw[very thick, latex-](t)--(cr'');
			\node(middle)[right =3cm of $(b1)!0.5!(bq'')$]{$\vdots$};
			\node(e2)[circle,fill,inner sep=0pt,minimum size=6pt, label = right:{$e_2$}][above = of middle]{};
			\node(e1)[circle,fill,inner sep=0pt,minimum size=6pt, label = right:{$e_1$}][above = of e2]{};
			\node(el)[circle,fill,inner sep=0pt,minimum size=6pt, label = right:{$e_m$}][below = of middle]{};
			\draw[very thick, -]($(e1)+(-0.85,0.425)$) rectangle node[above = 3cm]{$V_{BC}$}($(el)+(0.85,-0.425)$);
		\end{tikzpicture}}%
\end{center}
\caption{General representation of $H(G)$, as described in Section \ref{NWG}. (Internal arcs omitted.) Note that $m=|E_{BC}|$.} 
\label{fig:4}
\end{figure*}

\section{Minimum $\pmb{\T}$-transversal and Maximum Packing Equality}
\label{MTMPE}
 Before proving our main result, we recall two classical theorems that we will use. 
\vskip 6pt 
\noindent\textbf{\fontfamily{qhv}\selectfont{Menger's Theorem for Digraphs}} \cite{KM}
\hypertarget{MT}{}
\textit{Let $D$ be a finite directed graph and $s$ and $t$ be two vertices of $D$ with no arc from $s$ to $t$. Then the minimum size of an $(s,t)$-separator (that is, the minimum number of vertices, distinct from $s$ and $t$, whose removal breaks all directed paths from $s$ to $t$ in $D$) is equal to the maximum number of pairwise internally vertex-disjoint directed paths from $s$ to $t$ in $D$.}

\vskip 6pt
\noindent\textbf{\fontfamily{qhv}\selectfont{K\"onig's Line Colouring Theorem}} \cite{DK1}
\hypertarget{KT}{}
\textit{The minimum number of disjoint matchings that partition the edges of a bipartite graph equals the maximum degree of the graph.}
\vskip 12pt

\noindent We now prove Theorem~\ref{thm:1}: {\it If $G$ is a bilaterally-complete tripartite graph, then $\tau_{\triangle}(G)=\nu_{\triangle}(G)$.}

\begin{proof}
Let $G=(A,B,C;E)$ be a bilaterally-complete tripartite graph with $G[A\cup B]\cong K_{p,q}$ and $G[A\cup C]\cong K_{p,r}$, where $|A|=p$, $|B|=q$, $|C|=r$. We first show that every uniform $\T$-transversal of $G$ corresponds to an $(s,t)$-separator of $H(G)$ of the same size, and every $(s,t)$-separator of $H(G)$ corresponds to a (not necessarily uniform) $\T$-transversal of $G$ of the same size.

Let $E'=E'_{BC}\cup E_{AB'}\cup E_{AC'}$ be a uniform $\T$-transversal of $G$, where (as described in Section~\ref{SMTT}) $E'_{BC}\subseteq E_{BC}$, $B'\subseteq B$, $C'\subseteq C$, and $W=B'\cup C'$ is a vertex cover of $G[E_{BC}\setminus E'_{BC}]$. Then we construct an $(s,t)$-separator of size $|E'|$ in the network graph $H(G)$ as follows: Let $V'_{BC}\subseteq V_{BC}$ be the vertices in $H(G)$ corresponding to $E'_{BC}$, and for $1\leq i,j\leq p$ let $B'_i\subseteq B_i$ and $C'_j\subseteq C_j$ be the vertices in each copy of $B$ and $C$ corresponding to $B'$ and $C'$, respectively.

\textbf{Example:} The minimum $\T$-transversal shown in \hyperlink{G}{Figure 1} corresponds in \hyperlink{HG}{Figure 3} to the $(s,t)$-separator $\{2,2',38,9,9'\}$.
\vskip 6pt

\noindent \textbf{Claim:}  $\ds \V'=V'_{BC}\cup \bigcup_{i=1}^p B'_i \cup \bigcup_{j=1}^p C'_j$ is an $(s,t)$-separator of $H(G)$ of size $|E'|$.
\vskip 6pt

To see why this is true, note that every directed path from $s$ to $t$ in $H(G)$ contains one vertex in $B_i$, one vertex in $V_{BC}$, and one vertex in $C_j$, for some $i$ and $j$. Since $W=B'\cup C'$ covers all edges in $E_{BC}\setminus E'_{BC}$ (in $G$), every vertex of $H(G)$ in $V_{BC}\setminus V'_{BC}$ must be adjacent to at least one vertex in $B'_i\cup C'_j$, and therefore every path from $s$ to $t$ contains a vertex in $\V'$. We also have 
\[
|\V'|=|V'_{BC}|+\sum_{i=1}^p |B'_i|+\sum_{j=1}^p |C'_j|=|E'_{BC}|+p|B'|+p|C'|=|E'|.
\]

Conversely, let $\ds \V'=V'_{BC}\cup \bigcup_{i=1}^p B'_i \cup \bigcup_{j=1}^p C'_j$ be any $(s,t)$-separator in $H(G)$, with $V'_{BC}\subseteq V_{BC}$, $B'_i\subseteq B_i$, and $C'_j\subseteq C_j$. Every $s$--$t$ path in $\{s\}\cup B_i\cup V_{BC}\cup C_j\cup\{t\}$ must be cut by $\V'$, so in particular (when $i=j$), every vertex in $V_{BC}\setminus V'_{BC}$ must be adjacent to a vertex in $B'_i\cup C'_i$ for each $i=1,\ldots,p$.

By slightly abusing notation we will write $B'_i\subseteq B$ and $C'_i\subseteq C$ (in $G$) as well. (Note that we are not assuming that all of the $B'_i$ are the same set of vertices in $B$, for example.) Let $E'_{BC}$ be the edges in $ E_{BC}$ corresponding to the vertices in $V'_{BC}$. Then $W_i=B'_i\cup C'_i$ must cover all the edges in $ E_{BC}\setminus E'_{BC}$. So, labelling the vertices in $A$ as $a_1,\ldots, a_p$, we see that $\ds E'=E'_{BC}\cup \bigcup_{i=1}^p E_{a_iW_i}$ is a $\T$-transversal of $G$, and 
\[
|E'|=|E'_{BC}|+\sum_{i=1}^p|W_i|=|V'_{BC}|+\sum_{i=1}^p \left(|B'_i|+|C'_i|\right)=|\V'|.
\]
\vskip 12pt
Since there always exists a minimum $\T$-transversal of $G$ that is uniform, the above correspondences show that the smallest size of an $(s,t)$-separator in $H(G)$ equals $\tau_{\triangle}(G)$, the minimum $\T$-transversal size of $G$.
		 
Thus \hyperlink{MT}{Menger's Theorem} guarantees that there exists a set $P$ of $\tau_\triangle(G)$ pairwise internally vertex-disjoint $s$--$t$ paths in the network graph $H(G)$. Once again, each $s$--$t$ path in $H(G)$ contains one vertex from each of $B_i$, $V_{BC}$, and $C_j$ (for some $i,j$), which correspond to an edge in $E_{BC}$ and its endpoints. Let $F$ denote the subgraph of the bipartite graph $G[B\cup C]$ consisting of these edges and vertices corresponding to all paths in $P$. Because the paths in $P$ are internally vertex-disjoint, $F$ contains exactly $\tau_\triangle(G)$ edges. $F$ also has maximum degree at most $p$, because any vertex in $B$ or $C$ corresponds to only $p$ vertices in $H(G)$ and none of those vertices in $H(G)$ can appear more than once in the paths in $P$.
		 
		 Therefore by \hyperlink{KT}{K\"onig's Line Colouring Theorem}, the edges of $F$ can be partitioned into $p$ or fewer matchings, say $M_1,\ldots,M_k$ ($k\leq p$). We then join the edges in each matching $M_i$ to the vertex $a_i\in A$ for $1\leq i\leq k$, and the resulting set of triangles is clearly a packing in $G$. Since the number of edges in the graph $F$ was $\tau_\triangle(G)$, we have obtained a packing of size $\tau_\triangle(G)$ in $G$, which proves that $\tau_\triangle(G)=\nu_\triangle(G)$. 
\end{proof}

\section{Alternative Derivation}

The key step in our proof above is the existence of (up to) $p$ pairwise disjoint matchings in $G[B\cup C]$ that consist of a total of $\tau_\triangle(G)$ edges. This could actually be derived easily from Theorem~6 in Mao Cheng's 1983 paper \cite{CM} (in Chinese), as we show below. But that short paper gives only a brief proof based on Theorem~8 in the 1971 paper \cite{DS}, whose proof in turn relies on (essentially) Theorem~2.2 in the 1968 paper \cite{HP} (itself proved in a more general form from Menger's Theorem) and an unspecified result in \cite{OO}. It was due to the difficulty of tracing this convoluted logical route that we decided to present our own direct and self-contained proof above.

The statement of Mao Cheng's Theorem~6 is given below, with notation adjusted to correspond more closely with the present paper. Note that when $X$ is a set of vertices in a graph, $\partial(X)$ refers to the set of all edges in the graph that are incident on vertices in $X$. 

\vskip 6pt
\noindent\textbf{\fontfamily{qhv}\selectfont{Theorem 6}} {\fontfamily{qhv}\selectfont{(\cite{CM}, translated)}}
\hypertarget{CMT}{}
\label{CMT}
	\textit{Let $p,t\in\mathbb{N}^+$ and let $G[B\cup C]$ be a bipartite graph with parts $B$ and $C$. There exist $p$ pairwise edge-disjoint matchings in $G[B\cup C]$ such that the sum of their cardinalities equals $t$ if and only if for all $B''\subseteq B$ and for all $C''\subseteq C$,
		\begin{align*}
			|\partial(B'')\cap\partial(C'')|+p\left(|B|+|C|-|B''|-|C''|\right)\geq t.
		\end{align*}}

This result implies that $G[B\cup C]$ admits $p$ pairwise disjoint matchings with a total number of edges equal to the minimum value of 
\[\ds |\partial(B'')\cap\partial(C'')|+p\left(|B|+|C|-|B''|-|C''|\right)\]
 over all $B''\subseteq B$ and $C''\subseteq C$. We now show that this minimum value equals $\tau_{\triangle}(G)$ in the context of our tripartite graph $G=(A,B,C;E)$, which would allow us to complete the proof of Theorem~\ref{thm:1} by constructing a packing as before.

It was shown in Section~\ref{SMTT} that a (uniform) $\T$-transversal of $G$ corresponds to a set of edges $E'_{BC}\subseteq E_{BC}$ and subsets $B'\subseteq B$, $C'\subseteq C$ such that $W=B'\cup C'$ is a vertex cover of $G[E_{BC}\setminus E'_{BC}]$. If we let $B''=B\setminus B'$ and $C''=C\setminus C'$, then $E_{B''C''}=\partial(B'')\cap\partial(C'')\subseteq E'_{BC}$ is precisely the set of edges in $G[B\cup C]$ that are \emph{not} covered by $B'\cup C'$, so $E_{B''C''}\cup E_{AW}$ is a $\T$-transversal of $G$ of size
\[|E_{B''C''}|+p|W|\leq |E'_{BC}|+p|W|.\]
Conversely, for any $B''\subseteq B$ and $C''\subseteq C$, define $B'=B\setminus B''$, $C'=C\setminus C''$, and $W=B'\cup C'$. So $W$ is a vertex cover of $E_{BC}\setminus E_{B''C''}$, and $E_{B''C''}\cup E_{AW}$ is a uniform $\T$-transversal of $G$. 

Therefore the minimum value of 
\[ |\partial(B'')\cap\partial(C'')|+p\left(|B|+|C|-|B''|-|C''|\right)=|E_{B''C''}|+p|B'\cup C'|\] 
over all $B''\subseteq B$ and $C''\subseteq C$ is the minimum size of a uniform $\T$-transversal, which is $\tau_{\triangle}(G)$.\\

\noindent\textbf{\fontfamily{qhv}\selectfont{Acknowledgements }} We would like to thank You Rao for her translation of Mao Cheng's paper into English, and Donovan Hare for helpful comments on this manuscript.\\

\noindent\textbf{\fontfamily{qhv}\selectfont{Conflict of interest }} The authors declare that they have no conflict of interest.

\end{document}